# The Numerical Flow Iteration for the Vlasov–Poisson Equation

Matthias Kirchhart[*]    R. Paul Wilhelm[*]

**Abstract**

We present the numerical flow iteration (NuFI) for solving the Vlasov–Poisson equation. In a certain sense specified later herein, NuFI provides *infinite* resolution of the distribution function. NuFI exactly preserves positivity, *all* $L^p$-norms, charge, and entropy. Numerical experiments show no energy drift. NuFI is fast, requires several orders of magnitude less memory than conventional approaches, and can very efficiently be parallelised on GPU clusters. Low fidelity simulations provide good qualitative results for extended periods of time and can be computed on low-cost workstations.

## 1 Introduction

### 1.1 The Vlasov–Poisson Equation

Our problem of interest is a simplified model for the evolution of plasmas in their collisionless limit, as they occur in, for example, nuclear fusion devices. In dimensionless form it is given by the following system:

$$\partial_t f + v \cdot \nabla_x f - E \cdot \nabla_v f = 0, \tag{1}$$

$$E \coloneqq -\nabla_x \varphi, \tag{2}$$

$$-\Delta_x \varphi = \rho, \tag{3}$$

$$\rho(t,x) \coloneqq \bar\rho - \int_{\mathbb{R}^d} f(t,x,v)\,\mathrm{d}v. \tag{4}$$

Here, $f = f(t,x,v)$ is the electron distribution function, i.e., $f(t,x,v) \geq 0$ describes the probability density of electrons having velocity $v \in \mathbb{R}^d$ and location $x \in \mathbb{R}^d$ at time $t \in \mathbb{R}$. The Vlasov–Poisson equation for Plasmas is typically investigated for $d \in \{1,2,3\}$. In this work we will assume periodic boundary conditions in $x$, i.e., $\Omega \coloneqq [0,L]^d$ with $L > 0$, such that for all Cartesian basis vectors $e_i \in \mathbb{R}^d$ and all $k \in \mathbb{Z}$ one has $f(t, x + Lke_i, v) = f(t,x,v)$.

To obtain a well-posed problem, non-negative initial data of $f$ needs to be supplied, usually at time $t = 0$:

$$f(t=0, x, v) = f_0(x,v) \text{ for a given } x\text{-periodic } f_0 \in L^1(\Omega \times \mathbb{R}^d) \cap L^\infty(\mathbb{R}^d \times \mathbb{R}^d),\ f_0 \geq 0 \text{ a.\,e.} \tag{5}$$

In many computational benchmarks $f_0$ is a smooth function which is given as a relatively simple expression. Usually $f_0$ decays exponentially as $|v| \to \infty$, so in computational practice one pretends that one has $f(t,x,v) = 0$ whenever $|v| \geq v_{\max}$ for some user-defined parameter $v_{\max}$.

Equation (4) defines the charge density $\rho$; the parameter $\bar\rho$ stems from the assumption of a uniform ion-background and needs to be chosen such that overall neutrality is preserved, i.e., $\forall t \geq 0 : \int_\Omega \rho(t,x)\,\mathrm{d}x = 0$.

Neglecting collisions and the magnetic field, the Vlasov equation (1) then describes the evolution of $f$ under the influence of the self-consistent electric field $E = E(t,x)$, given in terms of the electric potential $\varphi = \varphi(t,x)$, which in turn is given as the solution of the Poisson equation (3).

---

[*]Applied and Computational Mathematics, RWTH Aachen University, Schinkelstraße 2, 52062 Aachen, Germany.
 E-Mail: kirchhart@acom.rwth-aachen.de and wilhelm@acom.rwth-aachen.de

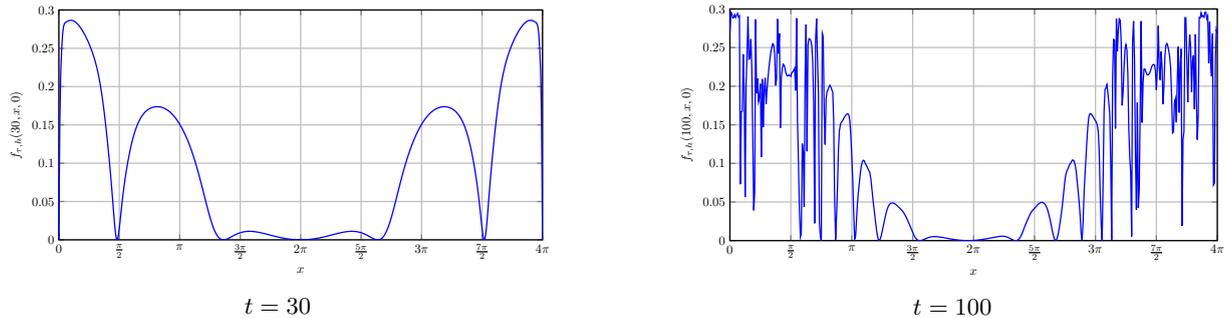

Figure 1: Cross-sections at $v = 0$ of approximations to distribution function $f(t,x,v)$ for the two stream instability benchmark in $d = 1$, see Section 4.2, at times $t = 30$ and $t = 100$. Computed using the numerical flow iteration as described below, using $N_x = 64$ and $N_v = 256$. The emergence of filaments which grow ever finer with time is clearly visible. Due to their finite resolution, all mesh-based approaches which directly discretise $f$ will eventually fail to reproduce these filaments accurately.

### 1.2 Challenges in Numerical Algorithms for the Vlasov–Poisson Equation

Superficially, numerically solving the Vlasov–Poisson equation might seem simple: assuming for the moment the electric field $E$ was known, the Vlasov equation (1) is a linear advection equation that can be solved using classical numerical schemes, e. g., finite differences with upstream discretisation. From such an approximation, the charge density $\rho$ could be computed directly; solving the Poisson equation (3) also is a classical problem for which there are many efficient algorithms available. The resulting approximation of the electric potential can then be fed back into the Vlasov solver to proceed by one time-step into the future.

While possible in theory, such a simple approach faces many problems in practice. Grid- and particle-based methods that directly discretise $f$ face the curse of dimensionality: the $2d$-dimensional $(x, v)$-space results in extremely large memory requirements. For such schemes, simulating the three-dimensional case ($d = 3$) usually is – if at all – only feasible on very large high-performance computers. At the same time, such schemes typically have very low arithmetic intensity, i. e., a low FLOP/byte count. Modern computer architectures struggle to deliver good performance for such schemes.

Even if these issues could be overcome, we want to particularly emphasise the difficulty due to so-called *filaments*, i. e., very fine details that develop in the solution over time. This is illustrated in Figure 1, which shows cross-sections of a numerical solution to a well-known benchmark problem at two different times $t$. The increasing filamentation in the solution manifests itself as oscillations with frequencies that rapidly increase over time $t$. As the filaments become too small for a given resolution, conventional discretisations often begin to violate several of the conservation laws discussed in Section 2.5. In this sense, these schemes then give unphysical results. For example, avoiding overshoots and negative values of $f$ near filaments is very difficult for many high-order schemes.

### 1.3 Overview

The rest of this article is structured as follows. Section 2 introduces the numerical flow iteration and discusses mathematical and practical issues. In Section 3 we elaborate on relationships with other methods found in the literature. Numerical experiments are discussed in Section 4. We conclude with an outlook to possible future extensions in Section 5.

## 2 The Numerical Flow Iteration

To understand NuFI, we will first need to introduce the exact flow associated with the Vlasov–Poisson equation. We then proceed by introducing the components of NuFI and also discuss several important conservation properties.

## 2.1 The *Exact* Flow and Solution to the Vlasov–Poisson Equation

Assume we knew the velocity $v \in \mathbb{R}^d$ and position $x \in \mathbb{R}^d$ of some imaginary 'particle' at time $s \in \mathbb{R}$. If we furthermore assume that the electric field $E$ was known for all times $t$, we can trace the state $\bigl(\hat{x}(t), \hat{v}(t)\bigr)$ of that imaginary particle both forward ($t > s$) and backward ($t < s$) in time by solving the following initial value problem:

$$\begin{aligned}
\frac{\mathrm{d}}{\mathrm{d}t}\hat{x}(t) &= \hat{v}(t), & \hat{x}(s) &= x, \\
\frac{\mathrm{d}}{\mathrm{d}t}\hat{v}(t) &= -E\bigl(t, \hat{x}(t)\bigr), & \hat{v}(s) &= v.
\end{aligned} \qquad (6)$$

Equation (6) is the definition of the characteristic lines of the Vlasov equation (1). Assuming $E$ is sufficiently smooth, as a consequence of the Picard–Lindelöf theorem, it is well-known that this system has a unique solution for any $(s, x, v) \in \mathbb{R} \times \mathbb{R}^d \times \mathbb{R}^d$. This leads us to the following definition:

**Definition 2.1 (Flow).** *The flow $\Phi_s^t(x, v)$ of the Vlasov–Poisson equation is defined as the unique solution of the initial value problem* (6):

$$\Phi_s^t(x, v) := \bigl(\hat{x}(t), \hat{v}(t)\bigr).$$

*We will sometimes use the term* exact flow *to distinguish $\Phi_s^t$ from approximative, numerical flows.*

In other words: $\Phi_s^t(x, v)$ tells us the position and velocity at some time $t$ of a particle with given state $(x, v)$ at time $s$. Note that system (6) is symplectic, and the following lemma is a well-known consequence that will be of great importance later on.

**Lemma 2.2 (Volume Preservation of the Exact Flow).** *The flow $\Phi_s^t : \mathbb{R}^d \times \mathbb{R}^d \to \mathbb{R}^d \times \mathbb{R}^d$ is a diffeomorphism that preserves volume in phase-space:*

$$\forall s, t \in \mathbb{R} : \ \forall (x, v) \in \mathbb{R}^d \times \mathbb{R}^d : \qquad \det \nabla_{(x,v)} \Phi_s^t(x, v) \equiv 1.$$

With help of the flow $\Phi_s^t$, the exact solution to the Vlasov–Poisson system takes a very simple form.

**Lemma 2.3 (Solution Formula).** *Let the non-negative distribution function $f$ be given at initial time $t = 0$, i.e., assume we were given $f_0(x, v) \geq 0$ such that $f(0, x, v) = f_0(x, v)$ for all $(x, v) \in \mathbb{R}^d \times \mathbb{R}^d$. Then the (exact) solution to the Vlasov–Poisson equation is given by:*

$$\forall (t, x, v) \in \mathbb{R} \times \mathbb{R}^d \times \mathbb{R}^d : \qquad f(t, x, v) = f_0\bigl(\Phi_t^0(x, v)\bigr).$$

## 2.2 The *Numerical* Flow and Solution to the Vlasov–Poisson Equation

In NuFI, the exact numerical flow $\Phi$ is replaced by a numerical approximation $\Psi$. For the moment we will still assume the electric field $E$ was known for all times $t$. For simplicity of the exposition, we will restrict ourselves to the case $\Phi_t^0 \approx \Psi_{n\tau}^0$, $n \in \mathbb{N}_0$, that is, we only consider discrete times $t = n\tau \geq 0$ with a user-defined time-step $\tau > 0$. However, extensions to arbitrary $t$ are certainly possible.

Due to the symplecticity of system (6), it is natural to approximate $\Phi_t^0$ using a symplectic one-step method for ordinary differential equations. In this work, we will use the well-known Störmer–Verlet method.[1] Noting that we are going backwards in time, in mathematical notation we have:

$$\Psi_{n\tau}^0 = \Psi_\tau^0 \circ \Psi_{2\tau}^\tau \circ \Psi_{3\tau}^{2\tau} \circ \Psi_{4\tau}^{3\tau} \circ \cdots \circ \Psi_{n\tau}^{(n-1)\tau}, \qquad (7)$$

where an individual step $\Psi_{k\tau}^{(k-1)\tau}(x, v)$, $k \in \{1, \ldots, n\}$, is given as:

$$\begin{aligned}
v_{k-\frac{1}{2}} &:= v + \tfrac{\tau}{2} E(k\tau, x), \\
x_{k-1} &:= x - \tau v_{k-\frac{1}{2}}, \\
v_{k-1} &:= v_{k-\frac{1}{2}} + \tfrac{\tau}{2} E\bigl((k-1)\tau, x_{k-1}\bigr), \\
\Psi_{k\tau}^{(k-1)\tau}(x, v) &:= (x_{k-1}, v_{k-1}).
\end{aligned} \qquad (8)$$

We believe that it is important to point out that $\Psi_{n\tau}^0$ can easily and efficiently be implemented on computer hardware, at a cost of just *one* evaluation of $E$ per time-step. In the following C++ snippet `vec` stands for a $d$-dimensional vector of floating point values.

```
    void numerical_flow( size_t n, double tau, vec &x, vec &v )
    {
        if ( n == 0 ) return;

        // Omit the following single line for Psi_tilda:
        v += (tau/2)*E(n*tau,x);

        while ( --n )
        {
            x -= tau * v;          // Inverse signs; we are
            v += tau * E(n*tau,x); // going backwards in time!
        }

        x -= tau*v;
        v += (tau/2)*E(0,x);
    }
```

This code also mentions the closely related function $\tilde{\Psi}^0_{n\tau}$, which is defined such that $\Psi^0_{n\tau}(x,v) = \tilde{\Psi}^0_{n\tau}\bigl(x, v + \frac{\tau}{2}E(n\tau,x)\bigr)$, which we will need later on.

The Störmer–Verlet method is second order accurate in time: $|\Phi^0_{n\tau}(x,v) - \Psi^0_{n\tau}(x,v)| = \mathcal{O}(\tau^2)$. Moreover, this method is symmetric, i.e., running the method backwards in time is the same as running the forward method on the time-reversed Vlasov–Poisson equation. With the numerical flow in place, it is natural to define approximations of the following shape:

$$f(n\tau,x,v) = f_0\bigl(\Phi^0_{n\tau}(x,v)\bigr) \approx f_0\bigl(\Psi^0_{n\tau}(x,v)\bigr) =: f_\tau(n\tau,x,v). \tag{9}$$

Several remarks are in order:

- Both $\Psi^0_{n\tau}$ and $\tilde{\Psi}^0_{n\tau}$ can be computed at arbitrary locations $(x,v) \in \mathbb{R}^d \times \mathbb{R}^d$. The Störmer–Verlet method only discretises in time, there is no limit to the $(x,v)$-resolution of the numerical flow. It is in this sense when we say that in NuFI the distribution function $f$ has *infinite resolution* in phase-space.

- It is not necessary to explicitly store $f$, only the initial data $f_0$ and the electric field $E$ need to be accessible.

- To evaluate $\Psi^0_{n\tau}$, the electric field $E$ is only required at the discrete time steps $i\tau$, $i = 0, \ldots, n$.

- To evaluate $\tilde{\Psi}^0_{n\tau}$, the electric field is only needed at previous times $i\tau$, $i = 0, \ldots, n-1$, i.e., the current $E(n\tau,\cdot)$ is not required.

In practice, the exact electric field $E(t,x)$ will be replaced with a numerical approximation $E_{\tau,h} \approx E$; the meaning of $h$ will be clarified later. The resulting fully discrete numerical flow $\Psi^0_{n\tau,h} \approx \Psi^0_{n\tau}$ will then also be distinguished with an additional index $h$. However, regardless of errors in $E_{\tau,h}$, just like the exact flow $\Phi^0_{n\tau}$, one of the remarkable key properties of symplectic integrators is that the numerical flows $\Psi^0_{n\tau}$ and respectively $\Psi^0_{n\tau,h}$ also preserve volume in phase-space. For completeness we repeat the simple proof.

**Lemma 2.4 (Volume Preservation of the Numerical Flow).** *Regardless of errors in the electric field $E_{\tau,h} \approx E$, the numerical flows $\Psi^0_{n\tau}$ and $\Psi^0_{n\tau,h}$ of the Störmer–Verlet method are diffeomorphisms which preserve volume in phase-space:*

$$\forall n \in \mathbb{N}_0 : \ \forall (x,v) \in \mathbb{R}^d \times \mathbb{R}^d : \quad \det \nabla_{(x,v)} \Psi^0_{n\tau}(x,v) \equiv \det \nabla_{(x,v)} \Psi^0_{n\tau,h}(x,v) \equiv 1.$$

*Proof.* We only describe $\Psi^0_{n\tau}$ as the proof for $\Psi^0_{n\tau,h}$ is identical. Because of the chain rule and the properties of the determinant one has:

$$\det \nabla_{(x,v)} \Psi^0_{n\tau} = \prod_{k=1}^n \det \nabla_{(x,v)} \Psi^{(k-1)\tau}_{k\tau}. \tag{10}$$

Similarly, a single step $\Psi^{(k-1)\tau}_{k\tau}$ is composed of the three sub-steps given in (8). The Jacobian matrix of, e.g., the first sub-step is given by:

$$\nabla_{(x,v)} \begin{pmatrix} x \\ v + \frac{\tau}{2} E(k\tau,x) \end{pmatrix} = \begin{pmatrix} I & 0 \\ \frac{\tau}{2} \nabla_x E(k\tau,x) & I \end{pmatrix}. \tag{11}$$

This is a triangular matrix with unit diagonal, so its determinant also is one. This is also the case for the other sub-steps, and thus $\det \nabla_{(x,v)} \Psi^0_{n\tau} \equiv 1$. $\square$

## 2.3 Structure of the Numerical Flow Iteration

So far we have assumed that the electric field $E$ was known exactly. In practice, however, this is of course not the case and it needs to be computed from the charge density $\rho$, which itself is defined in (4). Thus, only some approximation $E_{\tau,h} \approx E$ will be available, where $h$ denotes a discretisation parameter that will be described later.

At initial time $t = 0$ the charge density $\rho(t = 0, x)$ can be computed by integrating $f_0$ along $v$. This can, for example, be done by using numerical quadrature giving an approximation $\rho_{\tau,h}(0, x) \approx \rho(0, x)$. The approximate electric field $E_{\tau,h}(0, x)$ then results from the solution of the Poisson equation (3), which can itself be discretised with some reasonable numerical method.

From this point, we proceed iteratively in time. Thus, assume that we had already computed $\rho_{\tau,h}$, $\varphi_{\tau,h}$, and $E_{\tau,h}$ at times $i\tau$, $i = 0, \ldots, n-1$ for some $n \in \mathbb{N}$. We seek to compute $\rho_{\tau,h}(n\tau, x)$ using the numerical approximation $f(n\tau, x, v) \approx f_0\big(\Psi^0_{n\tau,h}(x, v)\big)$. However, $\Psi^0_{n\tau,h}$ cannot be evaluated, because $E_{\tau,h}(n\tau, x)$ is still unknown. Luckily, this circular dependency can be resolved by using $\tilde{\Psi}^0_{n\tau,h}$ and a simple change of variables:

$$\rho(n\tau, x) = \bar{\rho} - \int_{\mathbb{R}^d} f_0\big(\Phi^0_{n\tau}(x, v)\big)\,\mathrm{d}v$$
$$\overset{\Phi \leadsto \Psi_{\tau,h}}{\approx} \bar{\rho} - \int_{\mathbb{R}^d} f_0\big(\Psi^0_{n\tau,h}(x, v)\big)\,\mathrm{d}v = \bar{\rho} - \int_{\mathbb{R}^d} f_0\big(\tilde{\Psi}^0_{n\tau,h}(x, v + \tfrac{\tau}{2} E_{\tau,h}(n\tau, x))\big)\,\mathrm{d}v$$
$$\overset{\hat{v} := v + \tfrac{\tau}{2} E_{\tau,h}(n\tau, x)}{=} \bar{\rho} - \int_{\mathbb{R}^d} f_0\big(\tilde{\Psi}^0_{n\tau,h}(x, \hat{v})\big)\,\mathrm{d}\hat{v}. \quad (12)$$

In other words, for the computation of $\rho$, we can safely replace $\Psi^0_{n\tau,h}$ with $\tilde{\Psi}^0_{n\tau,h}$ without any additional error. Additionally, unlike $\Psi^0_{n\tau,h}$, its sibling $\tilde{\Psi}^0_{n\tau,h}$ *is computable*! By replacing the last integral with a suitable quadrature rule, we can thus compute $\rho_{\tau,h}(n\tau, x) \approx \rho(n\tau, x)$ and proceed. This is the numerical flow iteration, which can be summarised as follows:

1. Compute $\rho_{\tau,h}(t = 0, x)$ by integrating $f_0$ along $v$, using numerical quadrature.

2. Compute and store $\varphi_{\tau,h}(t = 0, x)$ using $\rho_{\tau,h}(t = 0, x)$.

3. For $n = 1, 2, 3, \ldots$:

    a) Compute an approximation $\rho_{\tau,h}(n\tau, x) \approx \rho(n\tau, x)$ using (12) and numerical quadrature.

    b) Compute and store $\varphi_{\tau,h}(n\tau, x)$ using $\rho_{\tau,h}(n\tau, x)$.

In the limit $h \to 0$, the last integral in (12) is evaluated without any error and the Poisson equation $-\Delta \varphi_{h,0} = \rho_{h,0}$ is also solved exactly.

**Definition 2.5 (Semi-discrete Approximation).** *The semi-discrete approximation $f_{\tau,0}$ is defined as*

$$f_{\tau,0}(n\tau, x, v) := f_0\big(\Psi^0_{n\tau,0}(x, v)\big),$$

*where at every time-step the last integral in (12) and the Poisson equation $-\Delta\varphi_{\tau,0} = \rho_{\tau,0}$ are solved exactly.*

Note that except for the initial time $t = 0$, we will still usually have $\Psi^0_{n\tau} \neq \Psi^0_{n\tau,0}$, $E \neq E_{\tau,0}$, $\varphi \neq \varphi_{\tau,0}$, and $\rho \neq \rho_{\tau,0}$ due to the errors introduced by the time discretisation. Consequently, we usually also have $f_{\tau,0} \neq f_\tau$ with $f_\tau$ from (9), which uses the exact field $E$.

## 2.4 Numerical Approximations of $\rho$ and $\varphi$

While the distribution function $f$ shows increasing filamentation over time, numerical evidence has shown that this is often not the case for the electric potential $\varphi$ and field $E$, and to a somewhat lesser extent also for the charge density $\rho$. An illustration of this is given in Figure 2.

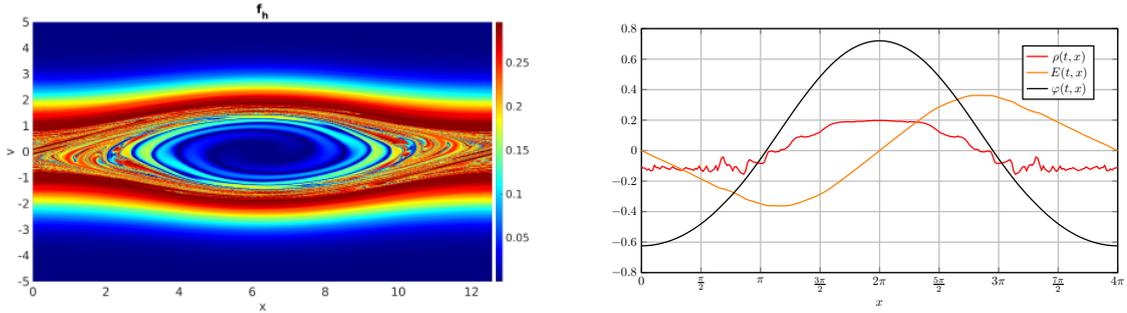

Figure 2: Even for very strongly filamented distributions $f$, the electric field $E$ and potential $\varphi$ often remain comparatively smooth. This example shows approximations of $f$ (left) as well as $\rho$, $E$, and $\varphi$ (right) for the two stream instability benchmark ($d = 1$) at $t = 100$. See Section 4.2 for more details.

For this reason, unlike $f$, the charge density $\rho$ and especially the electric potential $\varphi$ *can* efficiently be approximated on relatively coarse grids. In this work we sample the values of $\rho$ on a Cartesian grid of mesh-width $h_x = \frac{L}{N_x}$, for some $N_x \in \mathbb{N}$. For each grid node $x_i$, $i = 0, \ldots, N_x^d - 1$, we approximate $\rho(n\tau, x_i)$ using (12). The integral is approximated using the mid-point rule on a Cartesian grid reaching from $-v_{\max}$ to $v_{\max}$ in each of the $d$ components of $v$, having mesh width $h_v = \frac{2 v_{\max}}{N_v}$:

$$\rho(n\tau, x_i) \approx \rho_{\tau,h}(n\tau, x_i) := \bar{\rho} - h_x^d h_v^d \sum_{j=0}^{N_v^d - 1} f_0\bigl(\tilde{\Psi}^0_{n\tau,h}(x_i, v_j)\bigr). \tag{13}$$

Here, the $v_j$ are the mid-points of the cells of the Cartesian grid in $v$-direction. This requires the evaluation of $f_0 \circ \tilde{\Psi}^0_{n\tau,h}$ at $N_x^d N_v^d$ locations $(x_i, v_j)$ and is *by far* the computationally most expensive step of NuFI. Note however, that the evaluation of $f_0 \circ \tilde{\Psi}^0_{n\tau,h}$ is an embarrassingly parallel operation. The summation can thus be efficiently implemented on GPUs using so-called 'atomic additions'. It can also be parallelised to entire clusters of GPUs using simple reductions on partial results of $\rho_{\tau,h}$.

While in principle any quadrature rule could be used, the mid-point rule has the advantage of achieving exponential convergence for compactly supported smooth functions.[2] As mentioned in the introduction, in many cases $f_0$ is a smooth function that exponentially decays as $|v| \to \infty$ and can thus effectively be considered as compactly supported. Hence, before the occurrence of filamentation while $f$ still is smooth, we can expect very accurate results. At the same time, the mid-point rule has the advantage of being simple and stable: only positive terms appear in the sum.

We have thus defined $\rho_{\tau,h}$ at the grid-points $x_i$, between which we will interpolate. In this work we assume that $f$, and thus also $\rho$ and $\varphi$ are periodic in $x$. It thus makes sense to use trigonometric interpolation:

$$\rho_{\tau,h}(n\tau, x) = \sum_\alpha c_\alpha e^{i\alpha \frac{2\pi x}{L}}, \tag{14}$$

where $\alpha$ is a multi-index and the coefficients $c_\alpha$ can efficiently be computed from the point-values $\rho_{\tau,h}(n\tau, x_i)$ at the nodes $x_i$ using the fast Fourier transform (FFT). In this representation, solving the Poisson equation $-\Delta \varphi_{\tau,h} = \rho_{\tau,h}$ corresponds to a trivial scaling of the coefficients $c_\alpha$. The electric energy $\frac{1}{2} \|\nabla_x \varphi_{\tau,h}\|^2_{L^2(\Omega)}$ can also be computed easily using Parseval's identity.

While the spectral representation of $\varphi_{\tau,h}$ is essentially optimal from a mathematical point of view, it comes with a practical difficulty. When computing the numerical flow $\Psi^0_{n\tau,h}$ (or its relative $\tilde{\Psi}^0_{n\tau,h}$), the electric field and thus $-\nabla_x \varphi_{\tau,h}$ needs to be evaluated very often at *arbitrary* locations $x$. This evaluation is extremely expensive, as a summation over all modes $\alpha$ is necessary. For this reason, we follow a different approach.

Instead, we compute the inverse Fourier transform and thereby efficiently obtain the point values $\varphi_{\tau,h}(n\tau, x_i)$ at the grid nodes $x_i$, $i = 0, \ldots, N_x^d - 1$. Afterwards we compute the Cartesian, periodic, tensor-product spline interpolant of fourth order (piece-wise, coordinate-wise cubical polynomials of global coordinate-wise smoothness $C^2$). This representation is equally memory efficient, but a single evaluation of $-\nabla_x \varphi$ only takes a constant amount of time, independent of the resolution $h_x$, resulting in significant speed-ups.

**Definition 2.6 (Fully Discrete Approximation).** *The fully discrete approximation $f_{\tau,h}$ is defined as:*

$$f_{\tau,h}(n\tau, x, v) := f_0\big(\Psi^0_{n\tau,h}(x,v)\big),$$

*where $\Psi^0_{n\tau,h}$ is the numerical flow that is computed using the numerical approximation $E_{\tau,h} \approx E$ that results from taking the gradient of the spline interpolant of $\varphi_{\tau,h}$.*

## 2.5 Conservation Properties

What really makes NuFI stick out of the crowd, are its remarkable conservation properties. To avoid the technical details associated with the $x$-periodic setting, in this subsection we consider the whole-space case with $\bar{\rho} = 0$ instead, but note that these results also carry over to the periodic setting.

**Theorem 2.7 (Conserved Quantities).** *Let $F(n\tau, x, v)$ either denote the exact solution $F(n\tau, x, v) = f_0\big(\Phi^0_{n\tau}(x,v)\big)$, or a discrete approximation, $F(n\tau, x, v) = f_{\tau,h}(n\tau, x, v) = f_0\big(\Psi^0_{n\tau,h}(x,v)\big)$, $\tau > 0$, $h \geq 0$. Then, regardless of errors in the approximate electric field $E_{\tau,h} \approx E$ that is used in the computation of $\Psi^0_{n\tau,h}$, $F$ fulfils for all $n\tau$, $n \in \mathbb{N}_0$:*

- *the maximum principle: $0 \leq F(n\tau, x, v) \leq \|f_0\|_{L^\infty(\mathbb{R}^d \times \mathbb{R}^d)}$,*

- *conservation of all $L^p(\mathbb{R}^d \times \mathbb{R}^d)$ norms, $1 \leq p \leq \infty$:*

$$\|F(n\tau, \cdot, \cdot)\|_{L^p(\mathbb{R}^d \times \mathbb{R}^d)} = \|f_0\|_{L^p(\mathbb{R}^d \times \mathbb{R}^d)}, \tag{15}$$

- *conservation of kinetic entropy:*

$$\iint_{\mathbb{R}^d \times \mathbb{R}^d} F(n\tau, x, v) \ln\big(F(n\tau, x, v)\big) \, \mathrm{d}v\mathrm{d}x = \iint_{\mathbb{R}^d \times \mathbb{R}^d} f_0 \ln f_0 \, \mathrm{d}v\mathrm{d}x, \tag{16}$$

- *and, more generally, for any function $g : \mathbb{R} \to \mathbb{R}$ for which the following integrals are defined:*

$$\iint_{\mathbb{R}^d \times \mathbb{R}^d} g\big(F(n\tau, x, v)\big) \, \mathrm{d}v\mathrm{d}x = \iint_{\mathbb{R}^d \times \mathbb{R}^d} g\big(f_0(x,v)\big) \, \mathrm{d}v\mathrm{d}x. \tag{17}$$

*Proof.* For brevity, we will write $\Xi \in \{\Phi^0_{n\tau}, \Psi^0_{n\tau,h}\}$. The maximum principle directly follows from the fact that $F = f_0 \circ \Xi$ and $f_0 \geq 0$. Conservation of the $L^\infty$-norm follows because $\Xi$ is a diffeomorphism and thus $\Xi(\mathbb{R}^d \times \mathbb{R}^d) = \mathbb{R}^d \times \mathbb{R}^d$. Conservation of entropy ($g(x) = x \ln x$) and the other $L^p$-norms ($g(x) = x^p$) follows from the last statement. For the last statement, we transform the integral and use the volume preservation of $\Xi$ to obtain:

$$\iint_{\mathbb{R}^d \times \mathbb{R}^d} g\big(F(n\tau, x, v)\big) \, \mathrm{d}v\mathrm{d}x = \iint_{\mathbb{R}^d \times \mathbb{R}^d} g\big(f_0(\Xi(x,v))\big) \, \mathrm{d}v\mathrm{d}x =$$

$$\iint_{\mathbb{R}^d \times \mathbb{R}^d} g\big(f_0(x,v)\big) \underbrace{|\det \nabla_{(x,v)} \Xi^{-1}(x,v)|}_{\equiv 1} \, \mathrm{d}v\mathrm{d}x = \iint_{\mathbb{R}^d \times \mathbb{R}^d} g\big(f_0(x,v)\big) \, \mathrm{d}v\mathrm{d}x. \tag{18}$$

$\square$

The semi-discrete approximation $f_{\tau,0}$ additionally conserves momentum.

**Theorem 2.8 (Conservation of Momentum).** *Let $F(n\tau, x, v)$ either denote the exact solution $F = f_0\big(\Phi^0_{n\tau}(x,v)\big)$, or the semi-discrete approximation $F(n\tau, x, v) = f_{\tau,0}(n\tau, x, v)$ from Definition 2.5. Then $F$ fulfils the conservation of momentum:*

$$\forall n \in \mathbb{N}_0 : \quad \iint_{\mathbb{R}^d \times \mathbb{R}^d} vF(n\tau, x, v) \, \mathrm{d}v\mathrm{d}x = \iint_{\mathbb{R}^d \times \mathbb{R}^d} vf_0(x,v) \, \mathrm{d}v\mathrm{d}x. \tag{19}$$

*Proof.* The result for the exact solution is a classical matter; we thus only give the proof for $f_{\tau,0}$. For $n = 0$ the result is trivial as $f_{\tau,0}(0, x, v) \equiv f_0(x, v)$. Now assuming the result held for some $n \in \mathbb{N}_0$, we will show that it also holds for $k := n + 1$. Abbreviating $\Xi := \tilde{\Psi}^0_{k\tau,0}$, we obtain:

$$\iint_{\mathbb{R}^d \times \mathbb{R}^d} v f_{\tau,0}(k\tau, x, v) \, \mathrm{d}v \mathrm{d}x = \iint_{\mathbb{R}^d \times \mathbb{R}^d} v f_0\big(\Xi(x, v + \tfrac{\tau}{2} E_{\tau,0}(k\tau, x))\big) \, \mathrm{d}v \mathrm{d}x. \tag{20}$$

We now perform the same change of variables as in (12), however we note that due to the product $vf_{\tau,0}$ things become slightly more involved. Letting $\hat{v} = v + \frac{\tau}{2} E_{\tau,0}(k\tau, x)$, and noting that in this section we assume the whole-space case with $\bar{\rho} = 0$, we obtain:

$$\iint_{\mathbb{R}^d \times \mathbb{R}^d} \hat{v} f_0\big(\Xi(x, \hat{v})\big) \, \mathrm{d}\hat{v} \mathrm{d}x + \frac{\tau}{2} \int_{\mathbb{R}^d} E_{\tau,0}(k\tau, x) \rho_{\tau,0}(k\tau, x) \, \mathrm{d}x. \tag{21}$$

Note that $E_{\tau,0} = -\nabla \varphi_{\tau,0}$ and $\rho_{\tau,0} = -\Delta \varphi_{\tau,0}$. Moreover, for any function $\varphi$ that is smooth enough and decays sufficiently fast at infinity we obtain using integration by parts, interchanging the order of differentiation, and the fact that the Laplacian is self-adjoint:

$$\int_{\mathbb{R}^d} \nabla \varphi \Delta \varphi \, \mathrm{d}x = -\int_{\mathbb{R}^d} \varphi \nabla \Delta \varphi \, \mathrm{d}x = -\int_{\mathbb{R}^d} \varphi \Delta \nabla \varphi \, \mathrm{d}x = -\int_{\mathbb{R}^d} \Delta \varphi \nabla \varphi \, \mathrm{d}x. \tag{22}$$

Thus, this integral equals its own negative, and therefore the second integral in (21) is zero. As for the first integral, we perform another change of variables and let $\hat{x} = x - \tau \hat{v}$:

$$\iint_{\mathbb{R}^d \times \mathbb{R}^d} \hat{v} f_0\big(\Xi(x, \hat{v})\big) \, \mathrm{d}\hat{v} \mathrm{d}x = \iint_{\mathbb{R}^d \times \mathbb{R}^d} \hat{v} f_0\big(\Psi^0_{n\tau,0}(\hat{x}, \hat{v} + \tfrac{\tau}{2} E_{\tau,0}(n\tau, \hat{x}))\big) \, \mathrm{d}\hat{v} \mathrm{d}\hat{x}. \tag{23}$$

Thus, by performing a final change of variables on $\hat{v}$ and repeating the same arguments as above, we obtain in total:

$$\iint_{\mathbb{R}^d \times \mathbb{R}^d} v f_{\tau,0}\big((n+1)\tau, x, v\big) \, \mathrm{d}v \mathrm{d}x = \iint_{\mathbb{R}^d \times \mathbb{R}^d} v f_{\tau,0}\big(n\tau, x, v\big) \, \mathrm{d}v \mathrm{d}x. \tag{24}$$

Hence by induction, the result follows. $\square$

Finally, the exact solution also satisfies conservation of energy. Unfortunately, it is yet unclear to what extent this is the case for NuFI and we will point to our numerical experiments at the end of this article.

**Theorem 2.9 (Conservation of Energy).** *The exact solution $f(t, x, v) = f_0\big(\Phi^0_t(x, v)\big)$ and electric field $E$ satisfy the conservation of energy:*

$$\frac{\mathrm{d}}{\mathrm{d}t} \left( \frac{1}{2} \iint_{\mathbb{R}^d \times \mathbb{R}^d} v^2 f(t, x, v) \, \mathrm{d}v \mathrm{d}x + \frac{1}{2} \int_{\mathbb{R}^d} E(t, x)^2 \, \mathrm{d}x \right) = 0. \tag{25}$$

## 2.6 Complexity

### 2.6.1 Memory Complexity

Throughout this work we will assume that $f_0$ is given as a compact mathematical expression and can thus essentially be stored at zero cost. In this case, NuFI essentially only requires us to store $\varphi_{\tau,h}$ for each time-step.

It cannot be overstated that due to the missing $v$-dependence of $\varphi_{\tau,h}$, storing the spline coefficients of $\varphi$ on a grid requires several orders of magnitude less memory than storing $f$. Especially when $d = 3$, memory consumption can be reduced *millions of times*. Let $n \in \mathbb{N}$ denote the number of time-steps that should be computed. We then need to store $\varphi_{\tau,h}(k\tau, \cdot)$ for $k = 0, 1, \ldots, n$. This results in:

$$\text{Memory Requirement} = (n+1) N_x^d \times \text{size of one floating point value}. \tag{26}$$

| $d$ | $N_x = N_v$ | Memory of $\varphi_{\tau,h}$ per step | Memory for storing $f$ |
|---|---|---|---|
| 2 | 128 | 0.125 MiB | 2 048 MiB |
| 2 | 256 | 0.5 MiB | 32 768 MiB |
| 2 | 512 | 2 MiB | 524 288 MiB |
| 3 | 128 | 0.015 GiB | 32 768 GiB |
| 3 | 256 | 0.125 GiB | 2 097 152 GiB |
| 3 | 512 | 1 GiB | 134 217 728 GiB |

Table 1: Memory requirements for storing the electric potential $\varphi$ on a $d$-dimensional grid compared to directly storing the distribution function $f$ on a $2d$-dimensional grid for various example discretisations.

While in NuFI memory requirements do grow with the number of time-steps, the data required for an individual step is negligible compared to storing $f$ on a $2d$ dimensional grid, requiring $N_x^d N_v^d$ floating point values. Assuming 8 bytes per floating point value, the different memory requirements are illustrated in Table 1.

This means that, even for simulations with $d = 3$ and $N_x = 256$, $\varphi_{\tau,h}$ can be completely stored for several hundred time-steps in the memory of modern GPUs.

### 2.6.2 Computational Complexity

By far, the computationally most expensive operation in NuFI is the evaluation of $\rho_h(n\tau, x_i)$ using numerical quadrature. The cost of all other operations is negligible: quadrature is carried out on the $2d$ dimensional $(x, v)$-space while all other operations are carried out on the lower-dimensional $x$-space.

For simplicity, we will assume that $f_0$ can be evaluated efficiently at $\mathcal{O}(1)$ cost. In each time-step $k = 1, 2, \ldots, n$, the numerical flux $\Psi^0_{k\tau,h}$ needs to be evaluated at $N_x^d N_v^d$ quadrature nodes. Evaluating $\Psi^0_{k\tau,h}$ in turn requires $k + 1$ evaluations of the approximate electric field $E_{\tau,h}$; giving a total cost of approximately:

$$N_x^d N_v^d \sum_{k=1}^{n}(k+1) = \mathcal{O}\left(\frac{n^2}{2} N_x^d N_v^d\right) \tag{27}$$

evaluations of $E_{\tau,h}$. For this reason it is crucial to have efficient routines for the evaluation of $E_{\tau,h}$ available. In our approach using splines, the cost of a single evaluation of $E_{\tau,h}$ is independent of $n$, $N_x$ and $N_v$.

Thus, unlike conventional methods, the cost of NuFI grows quadratically with the number of time-steps. This might make the method look unattractive. However, NuFI has a much higher FLOP/byte ratio than approaches that directly store $f$, so it achieves much higher performance on modern computer systems. Additionally, also due to the low storage requirements, it easy to parallelise on clusters. Finally, no other method known to the authors conserves as many physical properties as NuFI does. It is maybe for this reason, as our numerical experiments will confirm, that NuFI only requires relatively coarse resolutions $N_x$, $N_v$ to achieve qualitatively good results.

## 3 Similarities and Differences to Related Literature and Methods

A key idea in NuFI is to approximate and follow the phase-flow $\Phi_s^t$ and use the method of characteristics to exploit exact conservation of values of $f$ along the phase-flow. In a sense this is a similar idea to Lagrangian schemes like Smooth Particle Hydrodynamics (SPH),[3],[4],[5] Particle-In-Cell (PIC),[5],[6],[7] and interpolation- and approximation-based particle methods.[8],[9] The key differences between particle methods like PIC and NuFI are the direct discretisation in phase-space and that one traces the characteristics forward in time instead of backwards.

In particle methods one tries to discretise the phase-space directly via sampling the initial distribution function $f_0$ at a set of points in the phase-space called 'particles'. The particles are then traced along their respective characteristics. For the non-linear case one then has to compute the electric field from these particles and their respective carried values of $f$. The different Lagrangian schemes basically differ in this aspect, i.e., how the electric field is computed from the samples. For the often-used PIC schemes the particles are mapped to a grid in $x$, thereby computing an approximation to $\rho$ from which $E$ can be obtained via a Poisson solver on the same grid.

| Hardware | CPU | GPU |
| --- | --- | --- |
| Workstation Laptop | Intel Xeon E-2276M @ 2.8 GHz, 6 cores, 2 threads per core | Nvidia Quadro T1000 Mobile, 4 GB DDR5-VRAM |
| Claix GPU Cluster | 2×Intel Xeon Platinum 8160 'Skylake' @ 2.1 GHz, 24 Cores each, per node | 2×Nvidia Volta V100-SXM2, 16 GB HBM2-RAM, per node |

Table 2: Overview of the hardware configurations used for the numerical experiments.

While in principle PIC has several desirable conservation properties, it struggles with excessive levels of noise caused by the fine filamentation and steep gradients in the solution $f$. This makes remeshing every few time-steps necessary.[10] Remeshing reduces the noise in the simulation but introduces a certain numerical diffusion similar to semi-Lagrangian methods.[11] Due to the steep gradients, methods that approximate the solution $f$ directly will suffer of overshoots in the numerical solution.[9],[12]

Related to Lagrangian methods and in particular PIC are the 'semi-Lagrangian' methods.[12],[13],[14],[15] Similar to purely Eulerian approaches[16],[17],[18] which discretise the distribution function directly and update the values on grid-nodes, semi-Lagrangian schemes also use a full grid-based discretisation of $f$ in the phase-space. However, to avoid too restrictive time-step constraints and to reduce numerical diffusion they trace the 'movement of grid-points' along the characteristics to use the values of the previous time-step approximation of $f$ for interpolation or approximation in the current time-step. The possibility of choosing large time-steps also exists in NuFI.

The use of direct discretisations in the whole phase-space via a grid is expensive as discussed in Section 2.6.1. Ways to overcome this are discussed in for instance the papers by Kormann on decomposition of the solution into a tensor train format and Einkemmer on comparing the discontinuous Galerkin with a spline-based semi-Lagrangian approach.[19],[20]

# 4 Numerical Experiments

In this section we describe the results of several numerical experiments. For this we used two different hardware configurations available to us: our local workstation laptop as well as the Claix GPU Cluster of RWTH Aachen University, see Table 2.

## 4.1 Weak Landau Damping ($d = 1$)

The first test case is commonly called *weak Landau Damping*. The initial condition is

$$f_0(x,v) := \frac{1}{\sqrt{2\pi}} e^{-\frac{v^2}{2}} \bigl(1 + \alpha \cos(kx)\bigr), \quad (x,v) \in [0, L] \times \mathbb{R} \tag{28}$$

with $k = 0.5$, $\alpha = 0.01$, $L = 4\pi$. The velocity space is cut at $v_{\max} = 10$. We use different spatial resolutions $N_x$ and $N_v$, but all simulations use the same time-step $\tau = \frac{1}{16}$.

The initial state (28) is a small perturbation to the Maxwellian distribution

$$f_M(v) = \frac{1}{\sqrt{2\pi}} e^{-\frac{v^2}{2}} \tag{29}$$

which is an equilibrium of the Vlasov–Poisson equation (1). It is known that this particular perturbation gets periodically damped to zero in a weak sense as $t \to \infty$. Thus, this benchmark consists of reproducing the correct damping rate for the electric field. The so-called Landau damping rate is $\gamma_E = 0.153359$, half of the corresponding damping rate for the electric energy $\gamma_{\text{energy}} = 0.306718$.[19]

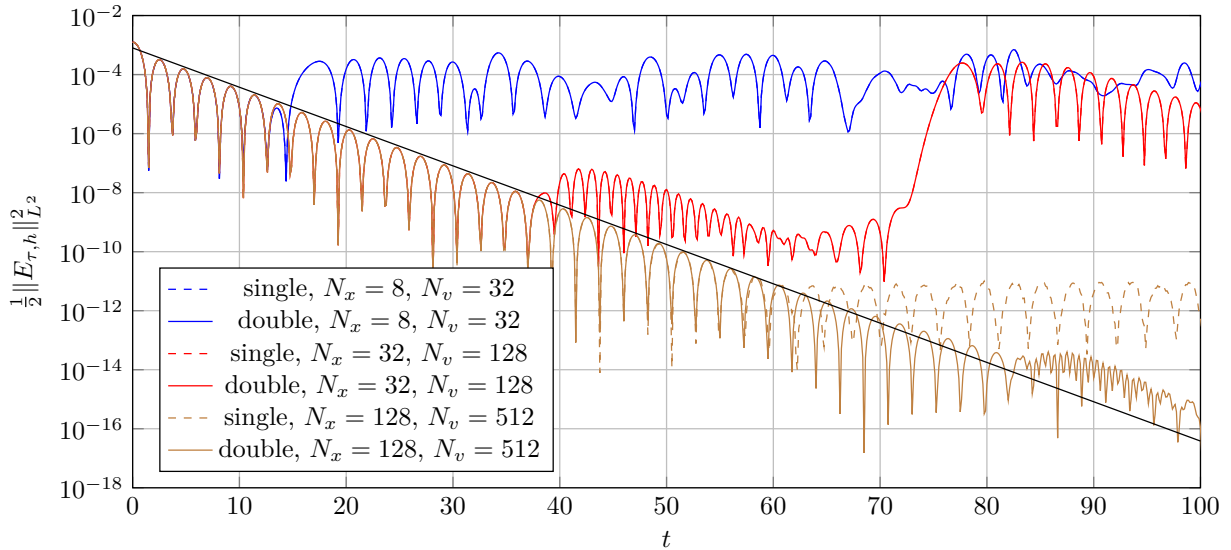

Figure 3: Electric energy for the *weak Landau damping* benchmark ($d = 1$) displayed for simulations run in single and double precision arithmetic. Except for the finest resolution, results for single and double precision are indistinguishable.

For $d = 1$ meaningful simulations can still be carried out on the workstation laptop. For such hardware computations in single precision are usually significantly faster and provide sufficient accuracy when $d = 1$. At the highest resolution, the computation of all 1 600 time-steps took 1.7 seconds in total when using single precision. As the laptop only uses a consumer grade GPU chip, the double precision computation was significantly slower and took 28.3 seconds.

The results are illustrated in Figure 3. For low and medium resolutions the results of single and double precision computations are indistinguishable. Only for high resolutions a different behaviour after $t \approx 60$ is observed: on the one hand, when computing in single precision no damping effect can be observed after $t \approx 62$ as the electric energy arrived at machine precision level with respect to the single precision, suggesting that Landau damping can at most be resolved up to the chosen floating point precision. On the other hand, the damping effect is still observable until $t \approx 80$ for double precision, however, with slightly slower rate. For single precision one does not observe recurrence phenomena: instead the electric energy oscillates with same amplitude after $t \approx 62$. The simulation with double precision shows a slight recurrence at $t \approx 83$.

The correct damping rate and oscillation frequency are reproduced by all simulations. The highest resolution simulation captures the correct damping rate $\gamma_{\text{energy}}$ until $t \approx 62$. Lower resolutions capture the correct damping rate only until their respective recurrence timings.

### 4.2 Two Stream Instability ($d = 1$)

Next we investigate the two stream instability benchmark. The initial condition is

$$f_0(x, v) := \frac{1}{\sqrt{2\pi}} e^{-\frac{v^2}{2}} v^2 \big(1 + \alpha \cos(kx)\big), \quad (x, v) \in [0, L] \times \mathbb{R} \qquad (30)$$

with $k = 0.5$, $\alpha = 0.01$, $L = 4\pi$. The velocity space is cut at $v_{\max} = 10$, time integration uses a time step of $\tau = \frac{1}{16}$. We again consider simulations in low ($N_x = 8$, $N_v = 32$), medium ($N_x = 32$, $N_v = 128$), and high resolution ($N_x = 128$, $N_v = 512$).

This benchmark simulates two colliding streams of electrons with opposing velocity which are both slightly perturbed from equilibrium. After an initial damping phase the streams start mixing, leading to an increasingly turbulent behaviour. In particular, the distribution function $f$ develops an increasing level of filamentation and a 'vortex' over time.

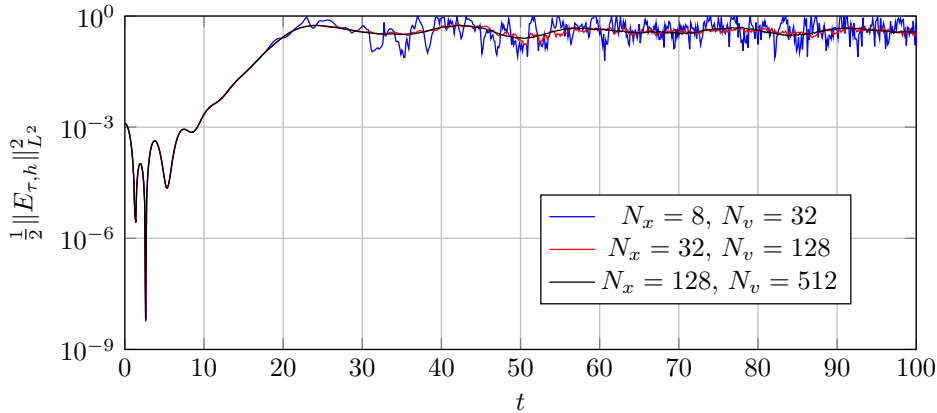

Figure 4: Electric energy for the *two stream instability* benchmark ($d = 1$) displayed for simulations, run in double precision and using time-step $\tau = \frac{1}{16}$.

Figure 4 shows the evolution of electric energy until $T = 100$. Until $t \approx 10$ damping effects on the two particle streams are still dominating the overall dynamic. Between $t \approx 10$ and $t \approx 25$ a vortex in the phase-space starts forming. After $t \approx 25$ the mixing process is dominant, where turbulence and filamentation in the solution make resolving $f$ increasingly complicated. Until $t \approx 25$ all three simulations are able to capture the underlying dynamics in the electric energy correctly. After $t \approx 25$ the low resolution simulation is still able to capture the dynamics qualitatively, however, with a significant error. After $t \approx 40$ the medium resolution simulation also starts to deviate from the high resolution one, albeit to a much lesser extent. The medium and high resolution simulation are able to resolve the slow oscillation in the electric energy after $t \approx 25$. The high resolution run is able to reproduce the electric energy until $t = 100$ without significant numerical noise in the plot.

Figure 6 shows the distribution function $f_{\tau,h}$ for several times $t$. At $t \approx 5$ the electric field damping seizes to be dominant. At this time the difference to the initial datum of $f$ is barely perceivable. Around $t \approx 25$ the forming of a vortex in phase-space starts setting in which comes with reaching the global maximum in electric energy. After $t \approx 25$ the vortex rotates periodically, creating an increasing number of filaments for later $t$. The low resolution is also able capture the overall dynamics correctly in a qualitative sense, however, one observes increasing distortions for the low resolution simulation. For $t = 100$ the vortex seems to move along the $x$-axis for the low resolution simulation, which is an incorrect dynamic.

The higher resolution simulation is able to correctly reproduce the dynamics until the end ($t = 100$). In addition, NuFI is able to reproduce the fine filamentation expected from the analytic solution. The pixel-artefacts which can be seen in Figure 6 are due to the finite resolution of the plot ($2\,048$ sampling points in each direction) and would disappear in plots of sufficiently high resolution.

Figure 5 shows the relative errors produced by simulation for several resolution settings in time and phase-space as well as different choice of floating point representation. The reference solution was computed using $N_x = 512$, $N_v = 2\,048$ and $\tau = \frac{1}{32}$.

When varying the phase-space resolution one observes that until $t \approx 20$ all errors stay at approximately the same level. Then the mixing process with accompanied with increasing filamentation in the solution sets in, which leads to an increase of error for all resolutions. After $t \approx 26$ errors seem to stabilise on a fixed level again; ranging from $10^{-3}$ to $10^{-5}$ for the relative $L^2$ error for the simulation in single precision and from $10^{-1}$ to $10^{-2}$ with respect to the $L^\infty$-norm. Similar error behaviour can be observed for double precision, suggesting that the chosen floating point precision makes no longer a difference when entering the turbulent phase of the simulation in this case. Note however, that for early times errors in double precision can be up to eight orders of magnitude smaller.

When varying the time-step-size and fixing the resolution in phase-space we observe an initial jump in error between $t = 0$ and $t \approx 2$. The error growth is much reduced for the smaller time steps. This implies that at this stage of the simulation the errors in the quadrature and the Poisson solver are close to machine precision $10^{-8}$, and that it is the time discretisation error which is dominating. After $t \approx 2$ the errors levels stabilise and do not significantly further increase until the onset of filamentation around $t \approx 26$.

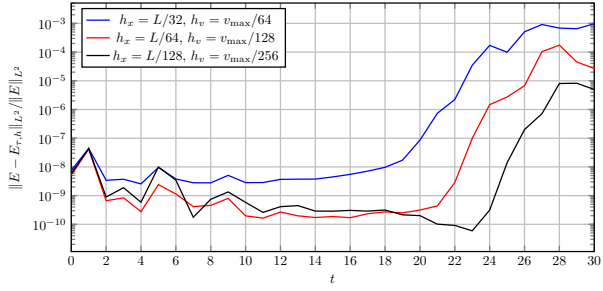

(a) $L^2$-error, single precision

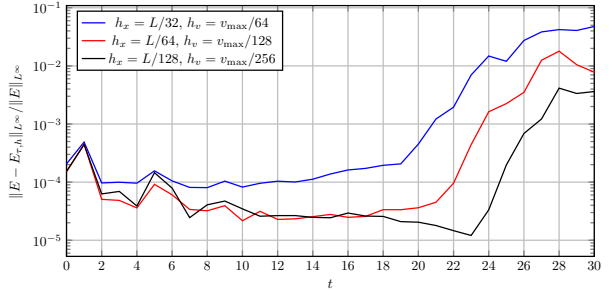

(b) $L^\infty$-error, single precision

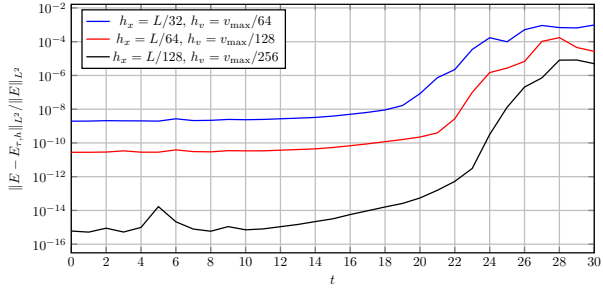

(c) $L^2$-error, double precision

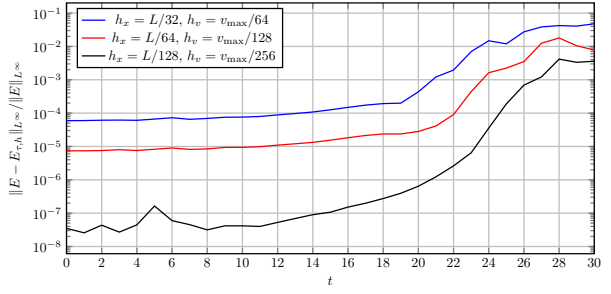

(d) $L^\infty$-error, double precision

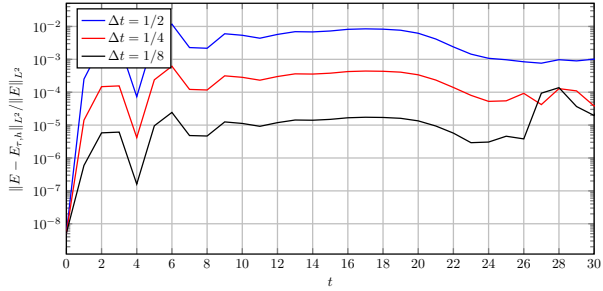

(e) $L^2$-error, single precision

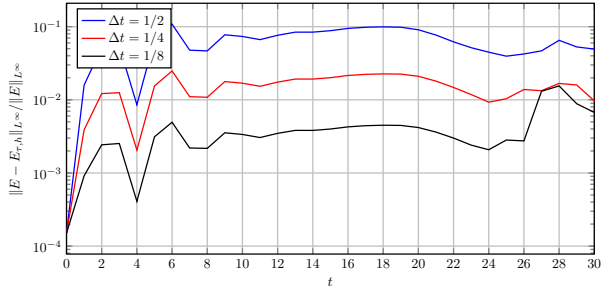

(f) $L^\infty$-error, single precision

Figure 5: Relative errors in electric field with respect to the $L^2-$ and $L^\infty-$norm for the *two stream instability* benchmark in $d=1$. With $\tau = \frac{1}{16}$ or $N_x = 64$, $N_v = 256$ fixed respectively. Reference solution computed using $N_x = 512$, $N_v = 2\,048$ and $\tau = \frac{1}{32}$.

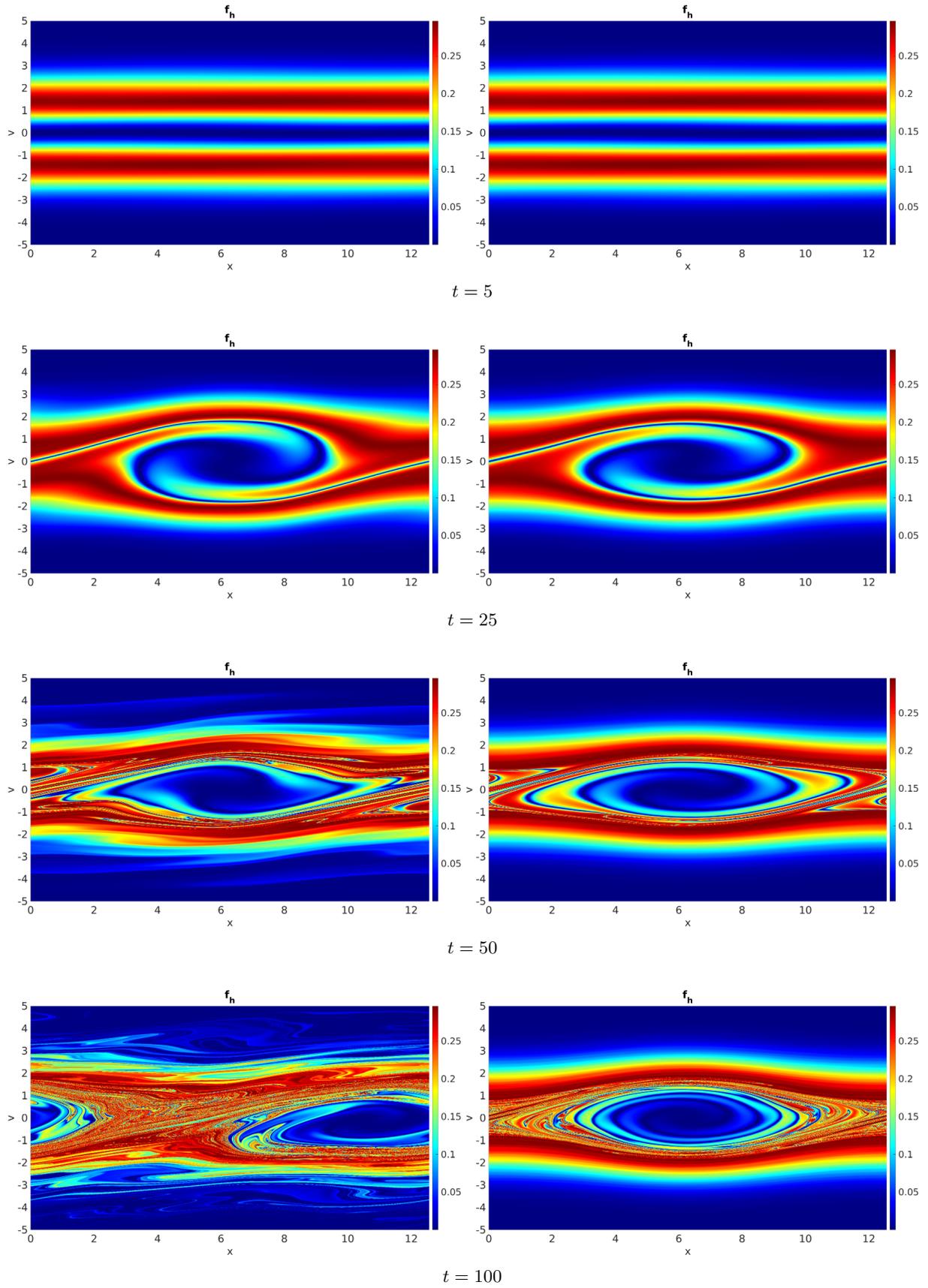

Figure 6: The distribution function $f_{\tau,h}(t,x,v)$ for the two stream instability benchmark. Left $N_x = 8$, $N_v = 32$ and right $N_x = 128$, $N_v = 512$. Both computed using $\tau = \frac{1}{16}$ as time-step. The plot resolution is 2 048 points in both $x$- and $v$-direction.

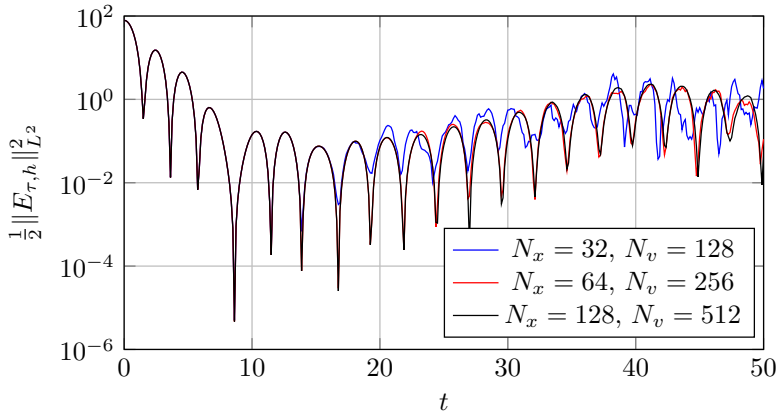

Figure 7: Electric energy for the *strong Landau damping* benchmark in $d = 2$, with time-step set to $\tau = \frac{1}{16}$.

## 4.3 Strong Landau Damping ($d = 2$)

For the two-dimensional case we consider the *strong Landau damping* or *non-linear Landau damping* benchmark. The initial condition is

$$f_0(x, y, u, v) := \frac{1}{2\pi} e^{-\frac{v^2 + u^2}{2}} \big(1 + \alpha(\cos(kx) + \cos(ky))\big), \quad (x, y, u, v) \in [0, L]^2 \times \mathbb{R}^2 \quad (31)$$

with $k = 0.5$, $\alpha = 0.5$, and $L = 4\pi$. The velocity space is cut at $v_{\max} = 10$, we chose $\tau = \frac{1}{16}$ as time-step. All computations are carried out using double precision arithmetic.

The initial datum is a strong perturbation of the Maxwellian distribution function. We thus expect a strong growth in the electric field strength. In this sense, this case is more comparable to the *two stream instability* than the *weak Landau damping* benchmark.

Figure 7 shows the evolution of electric energy over time. After a initial periodic damping period until $t \approx 10$ with electric energy ranging between $10^2$ and $10^{-1}$ the electric field gets excited again. The electric energy increases periodically until $t \approx 42$ after which it gets damped again.

While the high resolution simulation is able to resolve the dynamics in the electric energy correctly until $t = 50$, the low resolution simulation starts struggling after the excitement sets in, at $t \approx 20$ divergence from the high resolution solution becomes observable. In particular after $t \approx 30$ the low resolution simulation is also no longer able to capture the correct peaks and lows of the electric energy. However, the oscillating nature and overall dynamics are still reproduced qualitatively correctly.

In Figure 8 we plot relative errors with different resolutions in phase space. The reference solution was computed using $N_x = 256$, $N_v = 1\,024$ and $\tau = \frac{1}{32}$.

Initially the relative $L^2$-errors are again close to machine precision for the high resolution simulation, similar to the $d = 1$ case. For low resolution the simulation errors are around $10^{-8}$ with respect to the $L^2$-norm and $10^{-4}$ in the $L^\infty$-norm. However, in contrast to the $d = 1$ *two stream* benchmark the errors increase significantly faster. The errors increase continuously until $t \approx 20$ and level out at a range between $10^{-1}$ for low resolution in $L^2$- as well as $L^\infty$-norms and $10^{-4}$ in $L^2$- as well as $10^{-2}$ in $L^\infty$-norm for the high resolution simulation.

Figures 9 and 10 illustrate the electric field of high and low resolution simulations at various times $t$. There is no observable qualitative or quantitative difference between the two simulations in these plots, suggesting that even the low resolution is able able to correctly capture the dynamics in the electric field until at least $t = 50$.

## 4.4 Two Stream Instability ($d = 3$)

For the three-dimensional case we again consider the *two stream instability* benchmark. The initial condition is

$$f_0(x, y, z, u, v, w) := \frac{1}{2(2\pi)^{\frac{3}{2}}} \big(e^{-\frac{(v-v_0)^2}{2}} + e^{-\frac{(v+v_0)^2}{2}}\big) e^{-\frac{u^2+w^2}{2}} \big(1 + \alpha(\cos(kx) + \cos(ky) + \cos(kz))\big), \quad (32)$$

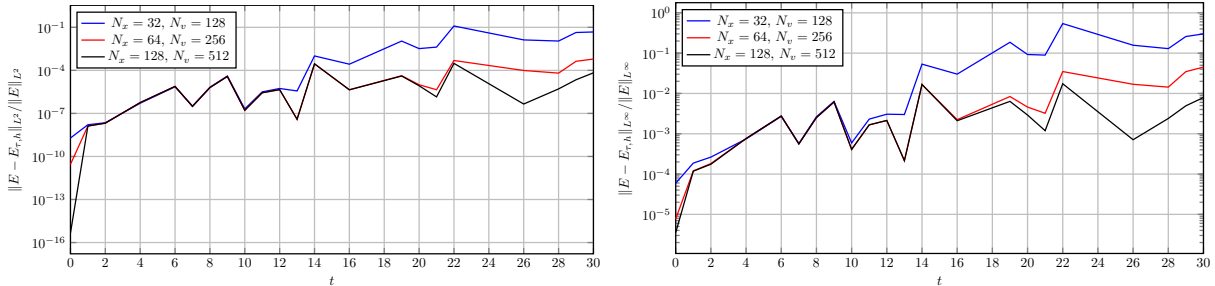

Figure 8: $L^2$- and $L^\infty$ errors with respect to a reference solution for the *strong Landau damping* benchmark in $d = 2$. Computed with time-step $\tau = \frac{1}{16}$. The reference solution was computed using $N_x = 256$, $N_v = 1\,024$.

with $(x, y, z, u, v, w) \in [0, L]^3 \times \mathbb{R}^3$ where $k = 0.2$, $\alpha = 10^{-3}$, $L = 4\pi$, and $v_0 = 2.4$. These parameters are the same as the ones chosen by Kormann.[19] The velocity space is cut at $v_{\max} = 10$ and we use $\tau = \frac{1}{16}$ as time-step. We carry out simulations at different resolutions; the finest uses $N_x = 32$, $N_v = 64$, and therefore *more than eight billion* quadrature nodes in phase-space. Using eight nodes of the Claix cluster with two GPUs each, the computation of all $n = 800$ time-steps took 3:58 hours. Again, only double precision arithmetic is considered.

Figure 11 shows the evolution of electric energy over time. In contrast to the $d = 1$ benchmark, the amplitude of the perturbation of the equilibrium is smaller, which leads to a delayed excitement of the electric field as well as delayed vortex-formation. The peak electric energy is reached at $t \approx 35$. The lower resolution simulation is not able to capture some of the peak-timings and after $t \approx 30$ exhibits a faster oscillation than the higher resolution simulation. However, even with low resolution of only 8 points in spatial and 16 points in velocity-direction the NuFI is still able to capture the overall dynamics in a qualitatively correct manner up to $t \approx 30$.

Figure 12 shows a representative cross-section of the electric field at $z = 0$. The electric field aligns along the $x$-axis with positive direction for approximately $y > 15$ and negative direction for $y < 15$.

### 4.5 On the Conservation of Energy

From Section 2.5 we know that all $L^p$-norms of $f_{\tau,h}$ as well as kinetic entropy are conserved exactly by NuFI. While *the conservation is exact*, the actual numerical values of the integrals in (15) and (16) can only be approximated by means of quadrature formulæ.

The same holds true for the total energy: its value can also only be approximated by means of quadrature. Additionally, however, we cannot expect an exact conservation of energy due to various discretisation errors. This even holds for the semi-discrete scheme where $h \to 0$ with a fixed time-step $\tau$. In particular, in many mechanical systems, symplectic time integrators usually do not conserve energy exactly; however the error in total energy remains bounded and does not grow over time. In other words: while these schemes usually do not conserve energy exactly, they are free of a systematic energy drift. Thus the question arises whether this is also the case in NuFI.

To at least partially address this question, we consider the two stream instability benchmark in $d = 1$, $\tau = \frac{1}{32}$. At time $t = 0$ the initial data $f_0$ and corresponding charge density $\rho$ are very smooth; thus the conserved quantities and the total energy can accurately be computed using numerical quadrature. The values obtained at $t = 0$ can thus serve as reference values. We then proceed as follows: in every time-step the integrals are evaluated using the same quadrature rule as in (13). The results are illustrated in Figure 13. For the entropy and $L^p$-norms the deviations from the value at $t = 0$ are quadrature errors only. The value for total energy contains both quadrature and discretisation errors. However, we observe that the relative errors of the approximate total energy are of the same magnitude as the quadrature errors for the other conserved quantities. We thus conclude that the error in total energy is smaller than the quadrature error; a systematic energy drift – if present – could only be detected by using more quadrature points. As this is the case for all tested resolutions, we conclude that an energy drift is unlikely.

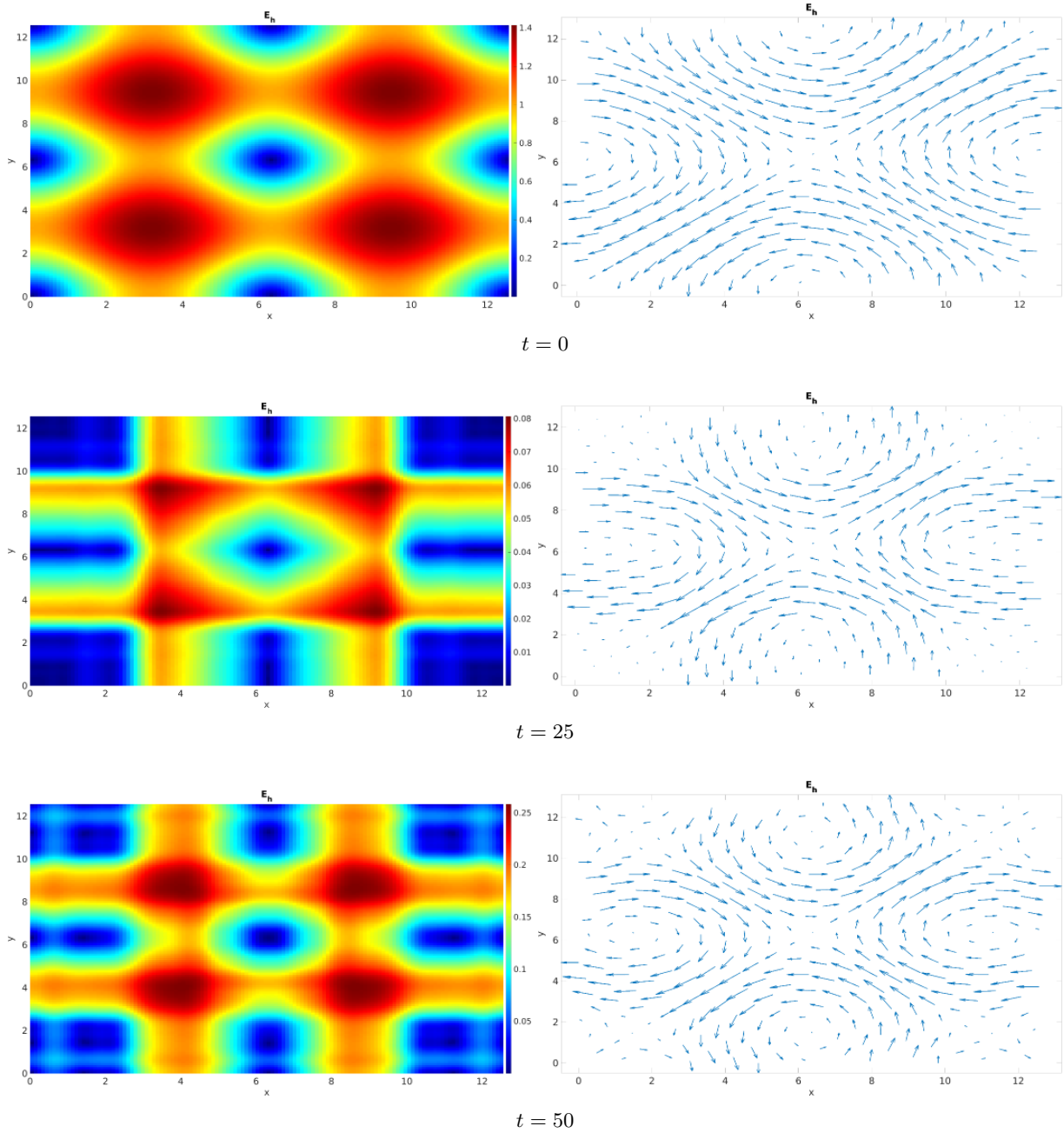

Figure 9: The electric field $E_{\tau,h}$ for the *strong Landau damping* benchmark in $d = 2$ at different times $t$. Left we display the magnitude $|E(t,x)|$ of the electric field, the right hand side illustrates the direction. The simulation was run with $\tau = \frac{1}{16}$, $N_x = 32$, $N_v = 128$.

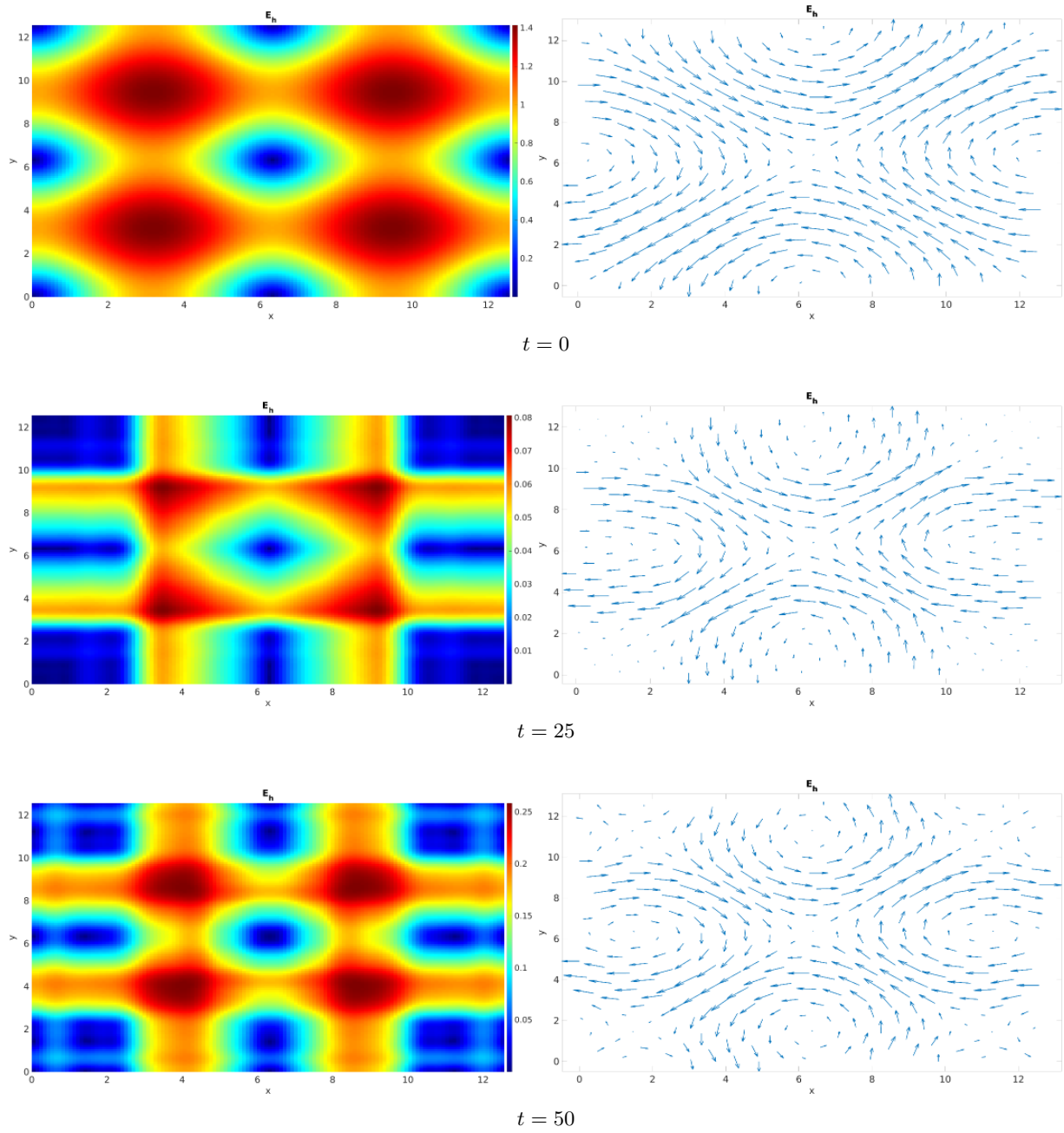

$t = 0$

$t = 25$

$t = 50$

Figure 10: The same as Figure 9, but this time with high resolution $N_x = 128$, $N_v = 512$. Both simulations qualitatively produce the same electric field.

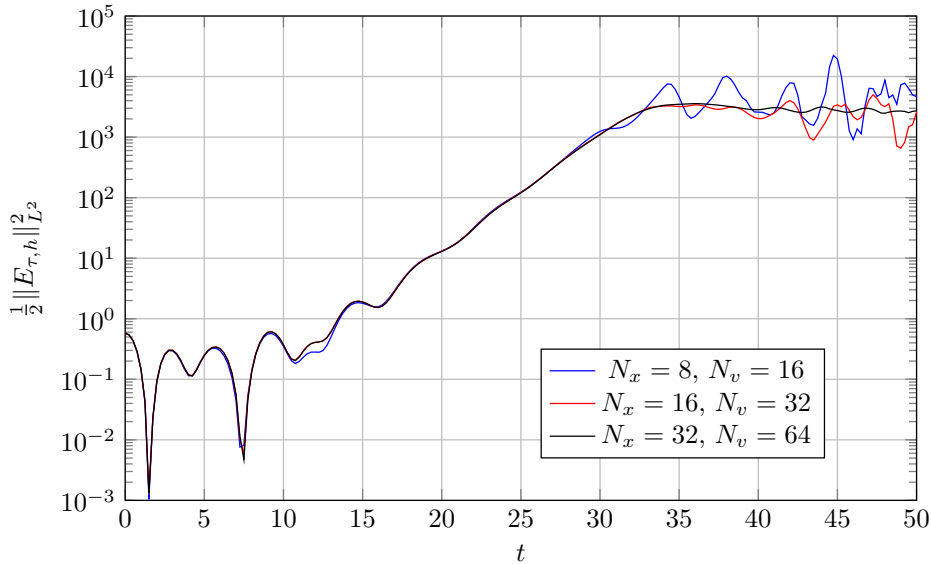

Figure 11: Electric energy for the *two stream instability* benchmark in $d=3$. Time-step $\tau = \frac{1}{16}$.

### 4.6 Computational Time and Scaling

As mentioned before, NuFI can be very efficiently implemented on GPU clusters. In this section, we thus look at the computational time and scaling behaviour of NuFI. The code for these computations is available from the authors per request. We consider the cases $d=1$ and $d=2$ only, as they allow larger numbers of samples: for $d=3$ the number of quadrature nodes grows too quickly and we would end up with graphs containing only a few computational results.

First we consider the computational time on a single node of the Claix cluster, see Table 2. For the case $d=1$ only one GPU and single precision was used, both GPUs and double precision was used for $d=2$. Starting with $N_x = 8$, $N_v = 32$, we measured the total computational time for $n=480$ time-steps of $\tau = \frac{1}{16}$, i.e., simulations were stopped at $t_{\text{end}} = 30$. For each subsequent simulation the values of $N_x$ and $N_v$ were doubled.

The computational times are displayed in Figure 14. For $d=1$ the computation time stays below $10^{-1}$ s until roughly one million quadrature nodes. Before this point, there are still spare execution units in the GPU thus computational time remains approximately constant. Afterwards we observe linear scaling until roughly one billion points. At this point there is again a upwards kink in the graph. We believe that this is due to caching effects: until this points the coefficients of $\varphi_{\tau,h}$ for a single time-step fit into the 128 KiB large cache of the GPU. For the $d=2$ we observe linear scaling starting from the second measurement point.

Next we consider a strong scaling experiment in $d=2$, with $N_x = 64$, $N_v = 256$, and $n=480$ time-steps, in double precision. Starting with a single node, we measure the computational time with an increasing number of computational nodes of the Claix cluster. The results are illustrated in Figure 15. We observe almost optimal strong scaling, i.e., when doubling the number of computational nodes the computation time is cut in half. After the measurement point with eight computational nodes there is slight kink in the plot which we believe is caused by the base calculations which take the same amount of time independent of the number of nodes like the Poisson solver, which uses a serial, non-parallelised implementation. We plan to carry out larger simulations on bigger machines, once corresponding computational resources become available to us.

## 5 Conclusion

The numerical flow iteration is a new method for numerically solving the Vlasov–Poisson equation with excellent conservation properties, low memory requirements, and high arithmetic intensity. It is ideally suited for modern ultra-parallel high-performance computing hardware and workstations alike. We believe that these properties will finally enable sufficiently resolved simulations of the Vlasov–Poisson equation for the most important case $d = 3$, already in 'Tier 2' data centres.

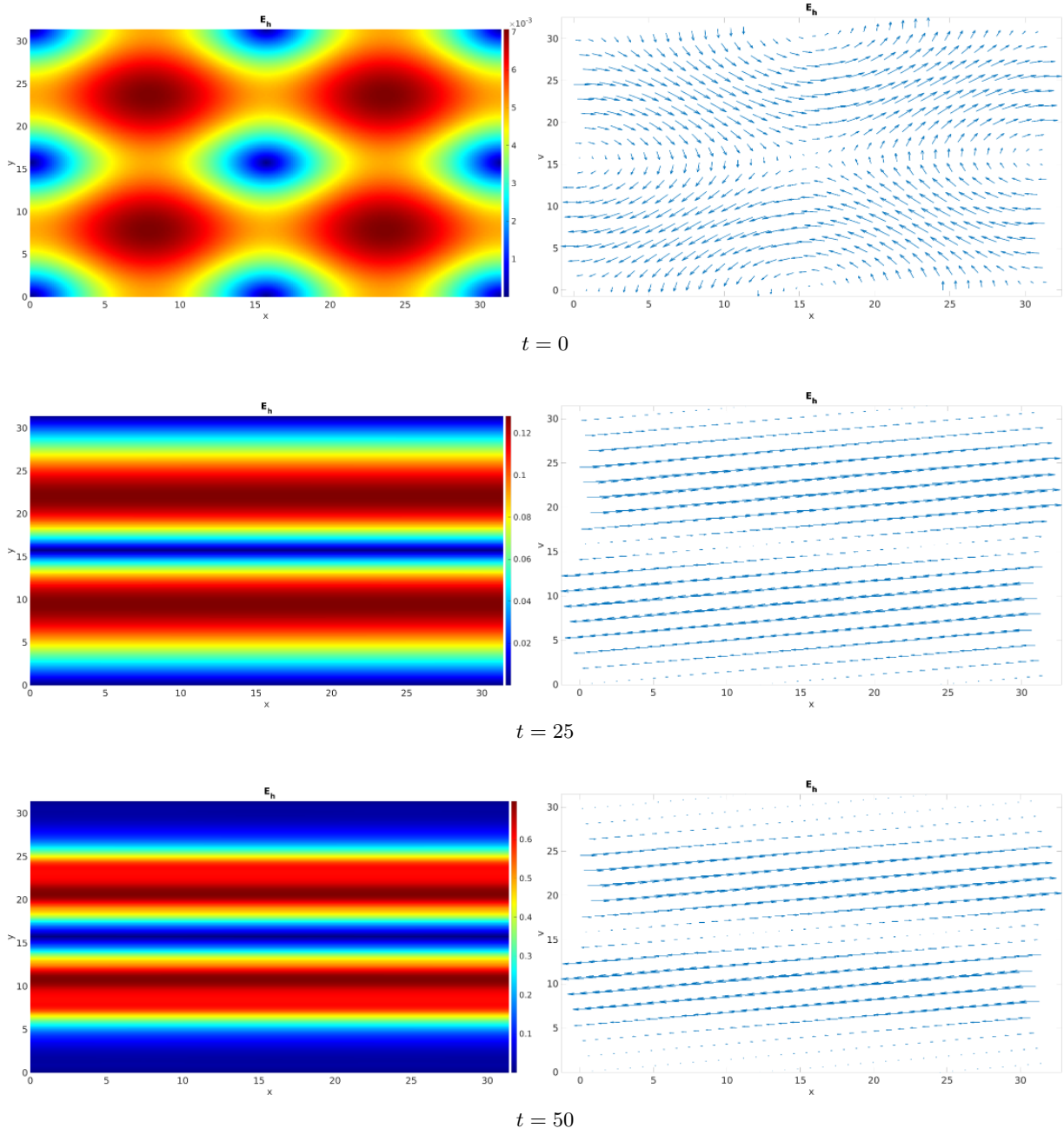

Figure 12: The electric field $E_{\tau,h}$ for the *two stream* benchmark in $d = 3$, at different times $t$ for $z = 0$. Left: magnitude $|E_{\tau,h}|^2$, right: directions in the $xy$-plane. The simulation was run with $\tau = \frac{1}{16}$, $N_x = 32$, $N_v = 64$, i.e., using 8 589 934 592 quadrature points.

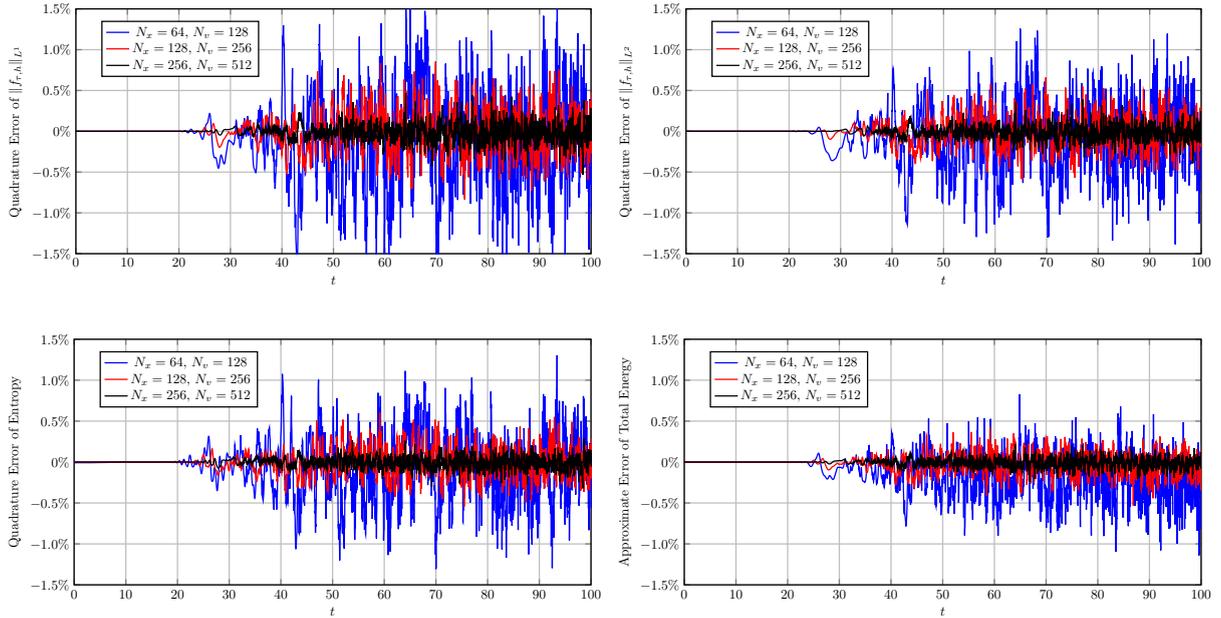

Figure 13: Quadrature errors when computing the $L^1$-norm, $L^2$-norm, and kinetic entropy of $f_{\tau,h}$ for the two stream instability benchmark in $d=1$ using $\tau = \frac{1}{32}$. *These quantities are conserved exactly*, their numerical value, however, can only be approximated using quadrature. The values are compared to the approximations at time $t=0$. Total energy cannot be expected to be conserved exactly, and cannot be evaluated exactly, either. However, its approximate error is of the same order of magnitude as the quadrature error of the exactly preserved quantities and no drift is visible.

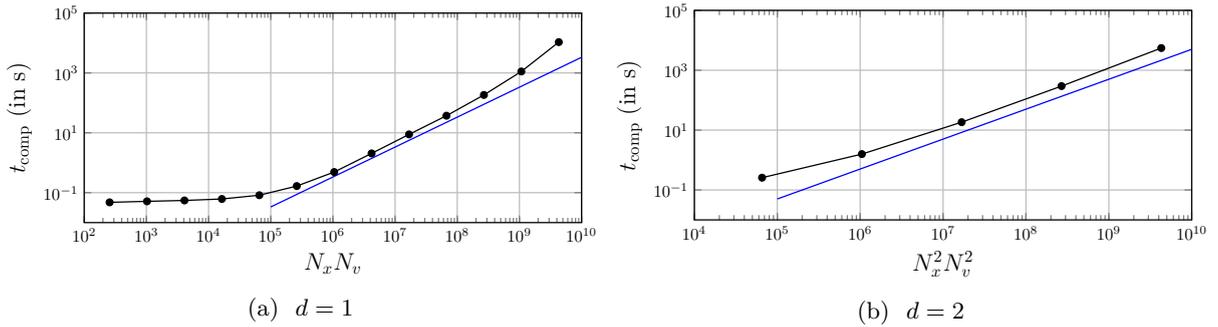

(a) $d=1$

(b) $d=2$

Figure 14: Computational times for $n=480$ time-steps and $N_x = 2^k$, $N_v = 2^{k+2}$, $k=3,4,5,\ldots$. Left: $d=1$, single precision, one GPU of the Claix cluster. Right: $d=2$, double precision, one node with two GPUs. Complexity grows linearly with the number of quadrature points $N_x^d N_v^d$. For $d=1$ caching effects become visible at $10^9$ quadrature points.

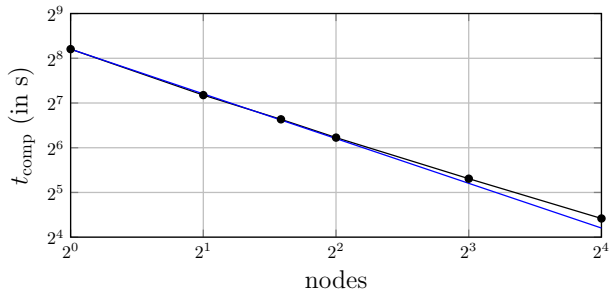

Figure 15: Strong scaling with the number of nodes for $d = 2$, $N_x = 64$, $N_v = 256$, $n = 480$ time-steps, and double precision, where both GPUs on each node of the Claix cluster are used. We observe almost optimal strong scaling; for 16 nodes parallel efficiency is only slightly reduced.

Several extensions to the scheme are possible. The approach can directly be applied to multi-species systems and can be extended to the full Maxwell equations.

In this work, for simplicity, we used entirely uniform discretisations. However, one can easily envision adaptive schemes: $\varphi_{\tau,h}$ can be stored using an adaptive space–time discretisation, adaptive numerical quadrature can be used, as well as adaptive time-stepping schemes. However, as usual, adaptive schemes make effective parallelisation significantly harder and require sophisticated load-balancing and synchronisation mechanisms.

Inhomogeneous right hand sides $g(t, x, v)$ of the Vlasov equation (1) can also easily be taken into account. In this case the solution formula from Lemma 2.3 becomes:

$$f(t,x,v) = f_0\big(\Phi_t^0(x,v)\big) + \int_0^t g\big(s, \Phi_t^s(x,v)\big)\,\mathrm{d}s. \tag{33}$$

The last integral can then be efficiently approximated using the numerical flow $\Psi_{\tau,h}$ and the trapezoidal rule. This way collision operators could also be incorporated into the simulation.

## Acknowledgments

Matthias Kirchhart started this work as a member of RWTH Aachen University, where he also received funding from the German research foundation (DFG), project number 432219818, 'Vortex Methods for Incompressible Flows'. The idea, analysis, and content were created there. He now works for Intel Deutschland GmbH, which kindly permitted the completion of this article.

Paul Wilhelm is recipient of a scholarship from the German national high performance computing organisation (NHR).

We would like to acknowledge the work of Jan Eifert, who helped in the implementation of the code.

Without these fundings and support, this work would not have been possible.